\documentclass{article}
\usepackage[utf8]{inputenc}
\usepackage{url}
\usepackage[margin=1.2in]{geometry}
\usepackage{amsmath,verbatim,amssymb,amsfonts,amscd,mathtools,xcolor}
\usepackage{algorithm,algpseudocode}
\usepackage{mathrsfs}
\usepackage{amsthm}
\usepackage{caption} 
\usepackage{subfigure}
\usepackage{lscape}
\usepackage[numbers]{natbib}
\usepackage{threeparttable}
\usepackage{booktabs}
\usepackage{comment}
\usepackage{adjustbox}
\usepackage{enumitem}

\newcommand{\be}{\begin{equation}}
\newcommand{\ee}{\end{equation}}
\newcommand{\ba}{\begin{aligned}}
\newcommand{\ea}{\end{aligned}}
\newcommand{\E}{\mathbb{E}}

\newcommand{\mtc}[1]{\mathcal{#1}}

\DeclarePairedDelimiterX\Set[2]{\lbrace}{\rbrace}%
 { #1 \,\delimsize| \,\mathopen{} #2 }
 
 \DeclareMathOperator*{\argmax}{arg\,max}
 
\theoremstyle{plain}

\title{A Model of Supply-Chain Decisions for Resource Sharing with an Application to Ventilator Allocation to Combat COVID-19}
\author{Sanjay Mehrotra, Hamed Rahimian, Masoud Barah, Fengqiao Luo, Karolina Schantz \\
\{mehrotra,\ hamed.rahimian,\ mbarah,\ fengqiaoluo2014,\ karolina.schantz\} \\
@northwestern.edu\\
    Industrial Engineering and Management Sciences Department \\
    Northwestern University, Evanston, IL 60208}

\date{April 01, 2020, Revised April 02, 2020}

\begin{document}
\maketitle

\begin{center}
    {\bf Abstract}
\end{center}
\noindent    
    We present a stochastic optimization model for allocating and sharing a critical resource in the case of a pandemic. The demand for different entities peaks at different times, and an initial inventory for a central agency is to be allocated. The entities (states) may share the critical resource with a different state under a risk-averse condition. The model is applied to study the allocation of ventilator inventory in the COVID-19 pandemic by FEMA to different US states. Findings suggest that if less than 60\% of  the ventilator inventory is available for non-COVID-19 patients, FEMA's stockpile of 20,000 ventilators (as of 03/23/2020) would be nearly adequate to meet the projected needs. However, when more than 75\% of the available ventilator inventory must be reserved for non-COVID-19 patients, various degrees of shortfall are expected. In an extreme case, where the demand is concentrated in the top-most quartile of the forecast confidence interval, the total shortfall over the planning horizon (till 05/31/20) is about 28,500 ventilator days, with a peak shortfall of 2,700 ventilators on 04/12/20. The results also suggest that in the worse-than-average to severe demand scenario cases, NY requires between 7,600-9,200 additional ventilators for COVID-19 patients during its peak demand. However, between 400 to 2,000 of these ventilators can be given to a different state after the peak demand in NY has subsided.
    
\vspace*{.2in}
\section{Introduction}
COVID-19 was first identified in Wuhan, China in December 2019 \cite{Huietal2020}. It has since become a global pandemic. As of March 31, 2020 the United States has overtaken China in the number of deaths due to the disease, with more than 3,900 deaths. Italy, which has 12,428 deaths, and Spain, which has 8,464, are the only two countries with higher death tolls. However, United States tops both of these countries in the current number of confirmed COVID-19 cases (189,035) \cite{JohnsHopkins2020}. Confirmed cases in the United States have more than doubled every three days in the time period since the first 100 cases were detected. This is even faster than the increases observed in Spain and Italy at the same point in the course of their epidemics \cite{OurWorldInDATA2020}. In Northern Italy, one of the global epicenters of the pandemic, COVID-19 has completely overwhelmed the healthcare system, forcing doctors into impossible decisions about which patients to save. Physicians on the front lines have shared accounts of how they must now weigh factors like age, comorbidities and probability of surviving prolonged intubation when deciding which patients with respiratory failure will receive mechanical ventilation \cite{Rosenbaum2020}. This experience is a warning of what awaits the United States. 

\subsection{A Resource Constrained Environment}

While approximately 80\% of COVID-19 cases are mild, the most severe cases of COVID-19 can result in respiratory failure, with approximately 5\% of patients requiring treatment in an intensive care unit (ICU) with mechanical ventilation  \cite{WooandMcGoogan2020}. Mechanical ventilation is used to save the lives of patients whose lungs are so damaged that they can no longer pump enough oxygen into the blood to sustain organ function. It provides more oxygen than can be delivered through a nasal cannula or face mask, allowing the patient’s lungs time to recover and fight off the infection. Physicians in Italy have indicated that critical COVID-19 patients often need to be intubated for a prolonged period of time (15-20 days) \cite{Rosenbaum2020}, further exacerbating ventilator scarcity. 

Limiting the death toll within the US depends on the ability to allocate sufficient numbers of ventilators to hard hit areas of the country before infections peak and ensuring that the inventory does not run out. Harder hit states (such as New York, Michigan and Louisiana) are now desperately trying to acquire additional ventilators in anticipation of significant shortages in the near future. Yet in the absence of a coordinated federal response, reports have emerged of states finding themselves forced to compete with each other in order to obtain ventilators from manufacturers \cite{CNBC2020}. According to New York's Governer Cuomo, the state has ordered 17,000 ventilators at the cost of \$25,000/ventilator, but is expected to receive only 2,500 over the next two weeks \cite{CuomoPressConference3-31-20}.  As of 03/31/2020, according to the US presidential news briefing, more than 8,100 ventilators have been allocated by FEMA around the nation. Of these, 400 ventilators have been allocated to Michigan, 300 to New Jersey, 150 to Louisiana, 50 to Connecticut, and 450 to Illinois, in addition to the 4,400 given to New York \cite{Trump033120}.

Going forward, the federal response to the COVID-19 pandemic will require centralized decision-making around how to equitably allocate, and reallocate, limited supplies of ventilators to states in need. Projections from the Institute for Health Metrics and Evaluation at the University of Washington, which assume that all states will institute strict social distancing practices and maintain them until after infections peak, show states will hit their peak demand at different time points throughout the months of April and May. Many states are predicted to experience a significant gap in ICU capacity, and similar, if not greater, gaps in ventilator capacity, with the time point at which needs will begin to exceed current capacity varying by state \cite{IHME2020}. 

\subsection{Our Contributions}
In response to the above problem, this paper presents a model for allocation and possible reallocation of ventilators that are available in the national stockpile. Importantly, computational results from the model also provide estimates of the shortfall of ventilators in each state under different future demand scenarios. 

This modeling framework can be used to develop master plans that will allocate part of the ventilator inventory here-and-now, while allocating and reallocating the available ventilators in the future. The modeling framework incorporates conditions under which part of the historically available ventilator inventory is used for non-COVID-19 patients, who also present themselves for treatment along with COVID-19 patients. Thus, only a fraction of the historical ventilator inventory is available to treat COVID-19 patients. The remaining demand needs are met by allocation and re-allocation of available ventilators from FEMA and availability of additional ventilators through planned production. The availability of inventory from a state for re-allocation incorporates a certain risk-aversion parameter. We present results while performing a what-if analysis under realistically generated demand scenarios using available ventilator demand data and ventilator availability data for different US states. An online planning tool is also developed and made available for use.  

\subsection{Organization}
This paper is organized as follows.  A review of the related literature is provided in Section~\ref{sec:Lit-Review}.
We present our resource allocation planning model,  and its re-formulation  in  Section~\ref{sec:Model}.   Section~\ref{sec:USCaseStudy}  presents  our  computational  results  under  different  mechanical ventilator demand scenarios for the COVID-19 pandemic in the US. This is followed by concluding remarks.

\section{Literature Review}\label{sec:Lit-Review}
A review on the role of operations research to ensure equity in global health is provided in \cite{bradley2017_OR-in-global-health}. The paper points out that poor availability of representative and high quality data, along with a lack of collaboration between operations research scientists, healthcare practitioners, and stakeholders are common challenges for effective operation research modeling in global health. A medical resource allocation problem in a disaster is considered in \cite{xiang2016_med-res-alloc-emerg}. Victims' deteriorating health conditions are modeled as a
Markov chain, and  the resources are allocated to optimize the total expected health recovery rate and reduce the total waiting time. Certain illustrative examples in a queuing network setting are also given in \cite{xiang2016_med-res-alloc-emerg}. The problem of scarce medical resource allocation after a natural disaster using a discrete event simulation approach is investigated in \cite{cao2012_princp-scarce-med-res-alloc-natural-disaster-sim-approach}. 
Specifically, the authors in \cite{cao2012_princp-scarce-med-res-alloc-natural-disaster-sim-approach} investigate four resource-rationing principles: first come-first served, random, most serious first, and least serious first. It is found that without ethical constraints, the least serious first principle exhibits the highest efficiency. However, a random selection provides a relatively fairer allocation of services and a better trade-off with ethical considerations. Resource allocation in an emergency department in a multi-objective and simulation-optimization framework is studied in \cite{feng2015_stoch-res-alloc-emerg-depart-multi-obj-sim-opt}. Simulation and queuing models for bed allocation are studied in \cite{vasilakis2001_sim-study-winter-bed-crisis,gorunescu2002_queue-model-bed-alloc}. 

%
%Ethical principles in fair allocation of scarce medical resources are studied in \cite{krutli2016-fair-alloc-scarce-med-res}. Results based on an online questionnaire are analyzed. Findings suggest that 'sickest first' and 'waiting list' receive the highest fairness endorsement by lay people and to some extent also by health professionals. Additional empirical studies on equitable and effective allocation of scarce medical resources can be found in  \cite{persad2009_alloc-scarce-med-interv,withanachchi2007_res-alloc-pub-hospital-effective}. 

 The problem of
determining the levels of contact tracing to control spread of infectious
disease using a simulation approach to a social network model is considered in \cite{armbruster2007_contact-tracing-contrl-infect-disease}. A linear programming model is used in 
investigating the allocation of HIV prevention funds across states \cite{earnshaw2007_lin-prog-alloc-HIV-funds-over-states}. This paper suggests that in the optimal allocation, the funds are not distributed in an equitable manner. A linear programming model
to derive an optimal allocation of healthcare resources in developing countries is studied in  \cite{flessa2000-lin-prog-alloc-med-res-developing-country}. Differential equation-based systems modeling approach is used in
\cite{araz2012_geo-priorit-distr-vaccines} 
to find a geographic and demographic dependent way of distributing pandemic influenza vaccines based on a case study of A/H1N1 pandemic. 

In a more recent COVID-19-related study,  the author \cite{kaplan2019_covid-19-prob-model} proposes a probability model to estimate the 
effectiveness of quarantine and isolation on controlling the spread of COVID-19. In  the context of ventilator allocation, a conceptual framework for allocating ventilators in a public emergency is proposed in \cite{ZazaEtal2016}. The problem of estimating mechanical ventilator demand in the United States during an influenza pandemic was considered in \cite{MeltzerEtAl2015}. In a high severity pandemic scenario, a need of 35,000 to 60,500 additional ventilators to avert  178,000 to 308,000 deaths was estimated. Robust models for emergency staff deployment in the event of a flu pandemic were studied in \cite{BienstockZenteno2012}. Specifically, the authors focused on managing critical staff levels during such an event,
with the goal of minimizing the impact of the pandemic. Effectiveness of the approach was demonstrated through experiments using realistic data.

A method for optimizing stockpiles of mechanical ventilators, which are critical for treating hospitalized influenza patients in respiratory failure, is introduced in \cite{HuangEtAl2017}. In a case-study, mild, moderate, and severe pandemic conditions are considered for the state of Texas. Optimal allocations prioritize local over central storage, even though the latter can be deployed adaptively, on the basis of real-time needs. Similar to this paper, the model in \cite{HuangEtAl2017} uses an expected shortfall of ventilators in the objective function, while also considering a second criteria of total cost of ventilator stockpiling. However, the model in \cite{HuangEtAl2017} does not consider distribution of ventilators over time. In the case of COVID-19,  the ventilator demand is expected to peak at different times in different states, as the demand for each state has different trajectories. Only forecasts are available on how the demand might evolve in the future. 

In this paper, we assume that the planning horizon is finite, and for simplicity we assume that reallocation decisions will be made at discrete times (days) $t_0, t_1, t_2, \dots$.  Under certain demand conditions, the ventilators may be in short supply to be able to meet the demand. Our model is formulated as a stochastic program, and for the purpose of this paper, we reformulate and solve the developed model in its extensive form. We refer the reader to 
 \cite{birge2011-stoch-prog-book,shapiro2014-stoch-prog-book} for a general description of this topic.
 
%In a multi-stage stochastic program, some random parameters are observed at different stages. At an intermediate stage, the decision maker should make the current optimal decision based on all decisions made at previous stages and observed random parameters. This multi-stage decision making feature exists in our ventilator allocation problem, where the random parameters are the demand, and the decision maker needs to decide how to distribute the ventilators to each state based on the demand observed in previous stages. Bender's decomposition \cite{birge2011-stoch-prog-book,bertsimas1997_introd-linear-prog}   and stochastic dual dynamic programming \cite{shapiro2011_analy-stoch-dual-DP} methods are developed to solve two-stage stochastic programs and multi-stage stochastic programs, respectively.
\section{A Model for Ventilator Allocation} 
\label{sec:Model}
In this section, we present a multi-period planning model to allocate ventilators to different regions, based on their needs, for the treatment of critical patients. We assume that the demand for ventilators at each planning period is stochastic. We further assume that there is a central agency that coordinates the ventilator (re)location decisions. The ventilators' (re)location is executed  at the beginning of a time period. Once these decisions are made and executed, the states can use their inventory to treat patients. Both the federal agency and the states have to decide whether to reserve their inventory in anticipation of future demand or share it with  other entities. 

Before presenting the formulation, we list the sets, parameters, and decision variables that are used in the model. 
\begin{itemize}
    \item Sets and indices 
        \begin{itemize}
            \item $\mtc{N}$: states (regions), indexed by $n \in \mtc{N}$,
            \item $\mtc{T}$: Planning periods, indexed by $t \in \mtc{T}$, 
        \end{itemize}
	\item Deterministic parameters
		\begin{itemize}
		    \item $T$: the total number of time periods, i.e, $\mtc{T}:=[T]$, where $[T]:=\{1, \ldots, T \}$,  
			\item $Y_n$: the initial inventory of ventilators in region $n \in \mtc{N}$ at time period $t=0$, 
			\item $I$: the initial inventory of ventilators in the central at the beginning of time period $t=1$, 
			\item $Q_t$: the number of ventilators produced during the time period $t-1$ that can be used at the beginning of time period $t$, for $t \ge 1$, 
			\item $\gamma_{n}$: the percentage of the initial  inventory of ventilators in region $n \in \mtc{N}$ that cannot be used to meet the demand for  patients at the critical level, 
			\item $\tau_{n}$: the percentage of the initial inventory of ventilators  in region $n \in \mtc{N}$ that the region is willing to share with other regions,  among those that can be used to care for patients at the critical level, 
			\item $\rho_n$: the risk-aversion of region $n \in \mtc{N}$ to send their idle ventilators to the central agency to be shared with other regions,
		\end{itemize}
	\item Stochastic parameter
	   \begin{itemize}
	       \item $\tilde{d}_{n,t}$: the number of patients in regions $n \in \mtc{N}$ at  the critical level that need a ventilator at the beginning of time period $t \in \mtc{T}$, 
	   \end{itemize}
	\item Decision variables 
		\begin{itemize}
			\item $x_{n,t}$: the number of ventilators relocated to region $n \in \mtc{N}$ from the central agency at the beginning of time period $t\in \mtc{T}$,
			\item $z_{n,t}$: the number of ventilators relocated to the central agency from region $n \in \mtc{T}$ at the beginning of time period $t \in \mtc{T}$,
			\item $y_{n,t}$: the number of ventilators at region $n \in \mtc{T}$ that can be used towards treating the patients at the  critical level at the end of time period $t \in \{0\} \cup \mtc{T}$, 
			\item $s_{t}$: the number of ventilators at  the central agency  at the end  of time period $t \in \{0\} \cup \mtc{T}$.
		\end{itemize}
\end{itemize}

The planning model to minimize the expected shortage of ventilators in order to treat patients at the critical level is formulated as follows: 

\begin{subequations}\label{opt:allocation}
%\addtolength{\jot}{1em}
\begin{align}
&\textrm{min}\;  \E\left[ \sum_{t \in \mtc{T}}  \sum_{n\in \mtc{N}}  (\tilde{d}_{n,t}-y_{n,t})^{+} \right]   & \label{eq: obj_shortage}\\
&\textrm{ s.t. } y_{n,t-1} + x_{n,t} - z_{n,t} = y_{n,t}, & \forall n\in   \mtc{N}, \; \forall t \in [T],  \label{eq: region_inv_balance}\\
& \qquad s_{t-1} +  Q_{t} + \sum_{n \in  \mtc{N}}z_{n,t}  - \sum_{n\in \mtc{N}} x_{n,t} = s_{t}, & \forall t\in[T], \label{eq: center_inv_balance} \\
& \qquad z_{n,t} \le \Big(y_{n,t} - (1-\tau_{n}) y_{n,0} - \rho_{n} \tilde{d}_{n,t} \Big)^{+},  & \forall n \in \mtc{N}, \; \forall t\in[T], \label{eq: region_safety_stock}\\
&\qquad \sum_{n \in  \mtc{N}}x_{n,t} \le s_{t-1} + Q_{t}+ \sum_{n \in  \mtc{N}}z_{n,t}, & \qquad \forall t \in [T],  \label{eq: center_relocation} \\
& \qquad y_{n,0}=(1-\gamma_{n}) Y_{n}, & \forall n\in\mtc{N}, \;  \label{eq: region_ini_inv}\\
& \qquad s_{0}= I, & \label{eq: center_ini_inv} \\
&\qquad x_{n,t}, \; z_{n,t} \ge 0,\;  & \forall n \in \mtc{N},\;\forall t\in[T], \label{eq: 1}\\
&\qquad y_{n,t} \ge 0, \;  & \forall n \in \mtc{N},\;\forall t\in\{0\} \cup [T], \label{eq: 2}\\
&\qquad s_{t}\ge 0, & \forall t\in\{0\} \cup [T] \label{eq: 3}. 
\end{align}
\end{subequations}
We will now explain  the model in detail. 
The objective function \eqref{eq: obj_shortage} denotes the expected total shortage of ventilators over all time periods $t \in \mtc{T}$ and all regions $n \in \mtc{N}$.  
Constraints \eqref{eq: region_inv_balance} and \eqref{eq: center_inv_balance} ensure the conservation  of ventilators for the regions and the central agency, respectively. Constraint \eqref{eq: region_safety_stock} enforces that a region is not sending out any ventilator  to the central agency if its in-hand inventory is lower   than its safety stock, where the safety stock is determined as $\rho_{n} \tilde{d}_{n,t}$, for $t \in [T]$ and $n \in \mtc{N}$. Constraint \eqref{eq: center_relocation} ensures that the total number of outgoing ventilators from the central agency to the regions  cannot be larger than the available inventory, after incorporating the newly produced ventilators and the incoming ones from other regions. Constraints \eqref{eq: region_ini_inv}  and \eqref{eq: center_ini_inv} set the initial inventory at the regions and central agency, respectively. The remaining constraints ensure the non-negativity of decision variables.

Note that the objective function and constraints \eqref{eq: region_safety_stock} are not linear. 
By introducing  an additional variable, the term $(\tilde{d}_{n,t}-y_{n,t})^{+}$ in the objective function, for $n \in \mtc{N}$ and $t \in \mtc{T}$, can be linearized as 
\begin{align*}
&  e_{n,t} \ge \tilde{d}_{n,t}-y_{n,t},  \\
& e_{n,t} \ge 0. 
\end{align*}
Furthermore, for each region $n \in \mtc{N}$ and time period $t \in \mtc{T}$, constraint \eqref{eq: region_safety_stock}  can be linearized as 
\begin{align*}
&  y_{n,t} - (1-\tau_{n}) y_{n,0}-\rho_{n} \tilde{d}_{n,t}\ge M(g_{n,t}-1),  \\
&   z_{n,t}\le y_{n,t} - (1-\tau_{n}) y_{n,0} - \rho_{n} d_{n,t} +M(1-g_{n,t}),\\
&  z_{n,t} \le M g_{n,t}, \\
& g_{n,t} \in \{0,1\},  
\end{align*}
where $M$ is a big number. 

As mentioned before, we assume that $\tilde{d}_{n,t}$, for $n \in \mtc{N}$ and $t \in \mtc{T}$, in model \eqref{opt:allocation} is a stochastic parameter. Let us suppose that $\tilde{d}_{n,t}$ has a finite support. This, in turns, implies that for each $t \in \mtc{T}$, the vectors $\mathbf{\tilde{d}}_{t}:=[\tilde{d}_{n,t}]_{n \in \mtc{N}}$   and $\mathbf{\tilde{D}}:=[\mathbf{\tilde{d}}_{t}]_{t \in \mtc{T}}$ have finite supports as well. We let $\Omega$ represent the finite support of $\mathbf{\tilde{D}}$, and use $\omega$ to denote an element of this set (i.e., a scenario).  Furthermore, suppose that $p^{\omega}$ represents the probability of scenario $\omega \in \Omega$, where $p^{\omega} \ge 0$,  and $\sum_{\omega \in \Omega} p^{\omega}=1 $. 

By incorporating the finiteness of the support of $\mathbf{\tilde{D}}$, a linearized reformulation of model  \eqref{opt:allocation}  can be written as  a mixed-binary program in the following extensive form:
\begin{subequations}\label{opt:allocation_linear}
\begin{align}
&\textrm{min}\; \sum_{\omega \in \Omega} p^{\omega} \left[\sum_{t \in \mtc{T}}  \sum_{n\in \mtc{N}}   e_{n,t}^{\omega} \right]   & \label{eq: obj_shortage_linear}\\
&\textrm{ s.t. } y_{n,t-1}^{\omega} + x_{n,t}^{\omega} - z_{n,t}^{\omega} = y_{n,t}^{\omega}, & \forall \omega \in \Omega, \; \forall n\in   \mtc{N}, \; \forall t \in [T],  \label{eq: region_inv_balance_omega}\\
& \qquad s_{t-1}^{\omega} +  Q_{t} + \sum_{n \in  \mtc{N}}z_{n,t}^{\omega}  - \sum_{n\in \mtc{N}} x_{n,t}^{\omega} = s_{t}^{\omega}, & \forall \omega \in \Omega, \;  \forall t\in[T], \label{eq: center_inv_balance_omega} \\
& \qquad y_{n,t}^{\omega} - (1-\tau_{n}) y_{n,0}^{\omega}-\rho_{n} \tilde{d}_{n,t}^{\omega}\ge M(g_{n,t}^{\omega}-1),  & \forall \omega \in \Omega, \; \forall n \in \mtc{N}, \; \forall t \in [T], \label{eq: lin_con_1}\\
& \qquad  z_{n,t}^{\omega}\le y_{n,t}^{\omega} - (1-\tau) y_{n,0}^{\omega} - \rho d_{n,t}^{\omega} +M(1-g_{n,t}^{\omega}), & \forall \omega \in \Omega, \;  \forall n \in \mtc{N}, \; \forall t \in [T], \label{eq: lin_con_2}\\
& \qquad z_{n,t}^{\omega} \le M g_{n,t}^{\omega}, & \forall \omega \in \Omega, \;  \forall n \in \mtc{N}, \; \forall t \in [T], \label{eq: lin_con_3}\\
&\qquad \sum_{n \in  \mtc{N}}x_{n,t}^{\omega} \le s_{t-1}^{\omega} + Q_{t}+ \sum_{n \in  \mtc{N}}z_{n,t}^{\omega}, &  \forall \omega \in \Omega, \; \forall t \in [T],  \label{eq: center_relocation_omega} \\
& \qquad y_{n,0}^{\omega}=(1-\gamma_{n}) Y_{n}, & \forall \omega \in \Omega, \;  \forall n\in\mtc{N}, \;  \label{eq: region_ini_inv_omega}\\
& \qquad s_{0}^{\omega}= I, & \forall \omega \in \Omega, \; \label{eq: center_ini_inv_omega} \\
& \qquad  e_{n,t}^{\omega} \ge d_{n,t}^{\omega}-y_{n,t}^{\omega}, &  \forall \omega \in \Omega, \; \forall n \in \mtc{N}, \; \forall t \in [T], \label{eq: lin_obj_omega}\\
&\qquad x_{n,t}^{\omega}, \; z_{n,t}^{\omega}, \; e_{n,t}^{\omega} \ge 0,\;  & \forall \omega \in \Omega, \;  \forall n \in \mtc{N},\;\forall t\in[T], \label{eq: 1_omega}\\
&\qquad y_{n,t}^{\omega} \ge 0, \;  & \forall \omega \in \Omega, \;  \forall n \in \mtc{N},\;\forall t\in\{0\} \cup [T], \label{eq: 2_omega}\\
&\qquad s_{t}^{\omega}\ge 0, & \forall \omega \in \Omega, \;  \forall t\in\{0\} \cup [T] \label{eq: 3_omega}, \\
& \qquad g_{n,t}^{\omega} \in \{0,1\}, & \forall \omega \in \Omega, \; \forall n \in \mtc{N}, \; \forall t \in [T], 
\end{align}
\end{subequations}
where $d_{n,t}^{\omega}$ denotes the number of patients at the critical level in regions $n \in \mtc{N}$ that need a ventilator at the beginning of time period $t \in \mtc{T}$ under scenario $\omega \in \Omega$. Note that all variables in model \eqref{opt:allocation_linear} have superscript $\omega$ to indicate their  dependence to scenario $\omega \in \Omega$. 

In our computational experiments in Section \ref{sec:USCaseStudy}, we used a commercial mixed-integer programming solver to obtain the results. 
Furthermore, we used $I+ \tau_{n} y_{n,0} +  \sum\limits{t^{\prime}  \le t} Q_{t}$ as a big-M  for $n \in \mtc{N}$ and $t \in \mtc{T}$. It is worth noting that \eqref{opt:allocation} (and \eqref{opt:allocation_linear} as well) considers multi-period decisions. In the model a decision maker will make  
decisions for the entire planning horizon using the information that is available at the beginning of planning. 

%this assumption. %The decision maker can solve a sequence of single-stage models at each time period in feature to obtain a heuristic decision for that time period.

\section{Ventilator Allocation Case Study: The US}
\label{sec:USCaseStudy}
The ventilator allocation model \eqref{opt:allocation_linear}, described in Section~\ref{sec:Model}, was implemented in Python 3.7. All computations were performed using GUROBI 9.1, on a  Linux Ubuntu environment, using  14 cores of a PC with $3.4$ GHz processor and $128$ GB of RAM. An hour time limit was given for all the runs.

\subsection{Ventilator Demand Data}
Since projected ventilator need is a key input for the model, it is important to use accurate estimates of the demand forecasts. The forecasts of ventilator needs generated by \cite{IHME2020} were used in our computational study. These forecasts were made available on 03/26/2020, and used the most recent epidemiological data and advanced modeling techniques. The available information closely tracks the real-time data \cite{IHME2020-OriginalSite}.  This COVID-19 needs forecast data was recently used in a recent presidential news brief \cite{Trump033120}. Although it is difficult to validate the ventilator need forecasts against actual hospital and state level operational data, as this information is not readily available, we find that this model's forecasts for deaths are quite accurate. For example, the model forecasted 217.9 deaths (CI: [176.95, 271.0]) on 03/29/20 for NY state. The number of reported deaths in the state on 03/29/30 were 237. Similarly, the model forecasted 262.2 deaths (CI: [206.9, 340]) on 03/30/20 against the actual deaths of 253 on that day. 

\subsection{Demand Scenario Generation}
We considered a seventy-day planning period, starting from March 23, 2020 and ending on May 31, 2020. 
We generated the random demands in ways that correspond to projected future demands under different mitigation effects. 
More precisely, we considered four different cases to generate random samples for the number of  ventilators needed to care for COVID-19 patients. These cases are listed below: 
\begin{itemize}[leftmargin=1.70cm]
    \item[Case I.] {\bf Average-I}: Each of the demand scenarios have equal probability and the distribution is uniform over the range of the CI provided in \cite{IHME2020},
    
    \item[Case II.] {\bf Average-II}:  The demand scenarios in the top 25\% of CI have 0.25 probability (equally distributed); and scenarios in the bottom 75\% have 0.75 probability,
    
    \item[Case III.] {\bf Worse than Average}: The demand scenarios in the top 25\% of CI have 0.50 probability; and  the scenarios in the bottom 75\% have 0.50 probability.
    
    \item[Case IV.] {\bf Severe}: The demand scenarios in the top 25\% of CI have 0.75 probability; and the  demand scenarios in the bottom 75\% have 0.25 probability,
    
    % \item[Case V.] {\bf Extreme }: The demand scenarios in top 25\% have 100\% probability, i.e., 0\% probability is assigned for the bottom 75\% demand scenarios.
\end{itemize}

We further discuss the demand generation procedure. A demand scenario contains the demand data for all days and states. In all Cases I--IV, we assumed that the forecast CI provided in \cite{IHME2020}, for each day and for each state, represents the support of the demand distribution. 

Case I and II are generated to develop average demand scenario representations that use the information provided in the CI given in \cite{IHME2020} in two different ways. In Case I, it is assumed that the mean is the median of the demand distribution (i.e., the right- and left-tail of the demand distribution have 0.5 probability). We randomly generated a  number to indicate which tail to sample from, where both tails have the same 0.5 probability of being chosen. Once the tail is determined, we divided  the tail into $50$ equally-distanced partitions, and chose a random partition to uniformly sample from.  We repeated this process for all days and states. We sampled from the same  tail and partition for all days and states, although the range from which we sample depends on the CI. In this case, all scenarios are equally likely.

In Case II, we randomly generated a  number to indicate which tail to sample from, where the top 25\% of the CI (i.e., the right tail) has a 0.25 probability and the bottom 75\%  (i.e., the left tail) has a 0.75 probability of being chosen. If the right tail is chosen, we set the weight of the scenario to 0.25, and we set it to 0.75 otherwise. The rest of the procedure is similar to Case I. In order to determine the probability of scenarios, we normalized the weights. 
Demand scenarios in Cases III-V are generated in the same fashion as in Case II, where the only difference is in the probability of which tail to choose from, which is determined by the sampling scheme described in the definition of the case. 

\begin{figure}
	\centering
	\footnotesize
	% \begin{minipage}{.7\linewidth}
	\begin{tabular}{cc}
		\includegraphics[width=65mm]{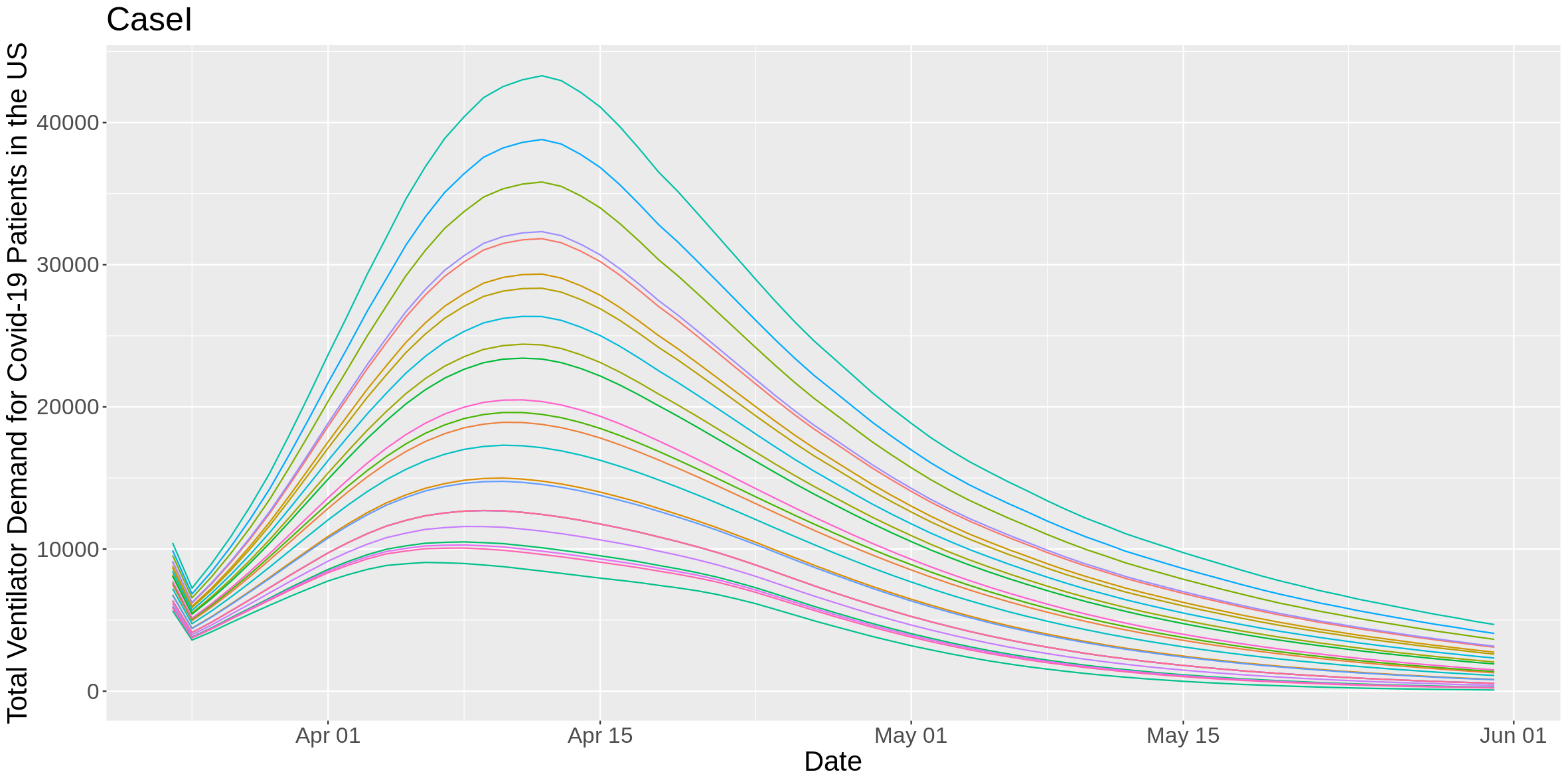} &   \includegraphics[width=65mm]{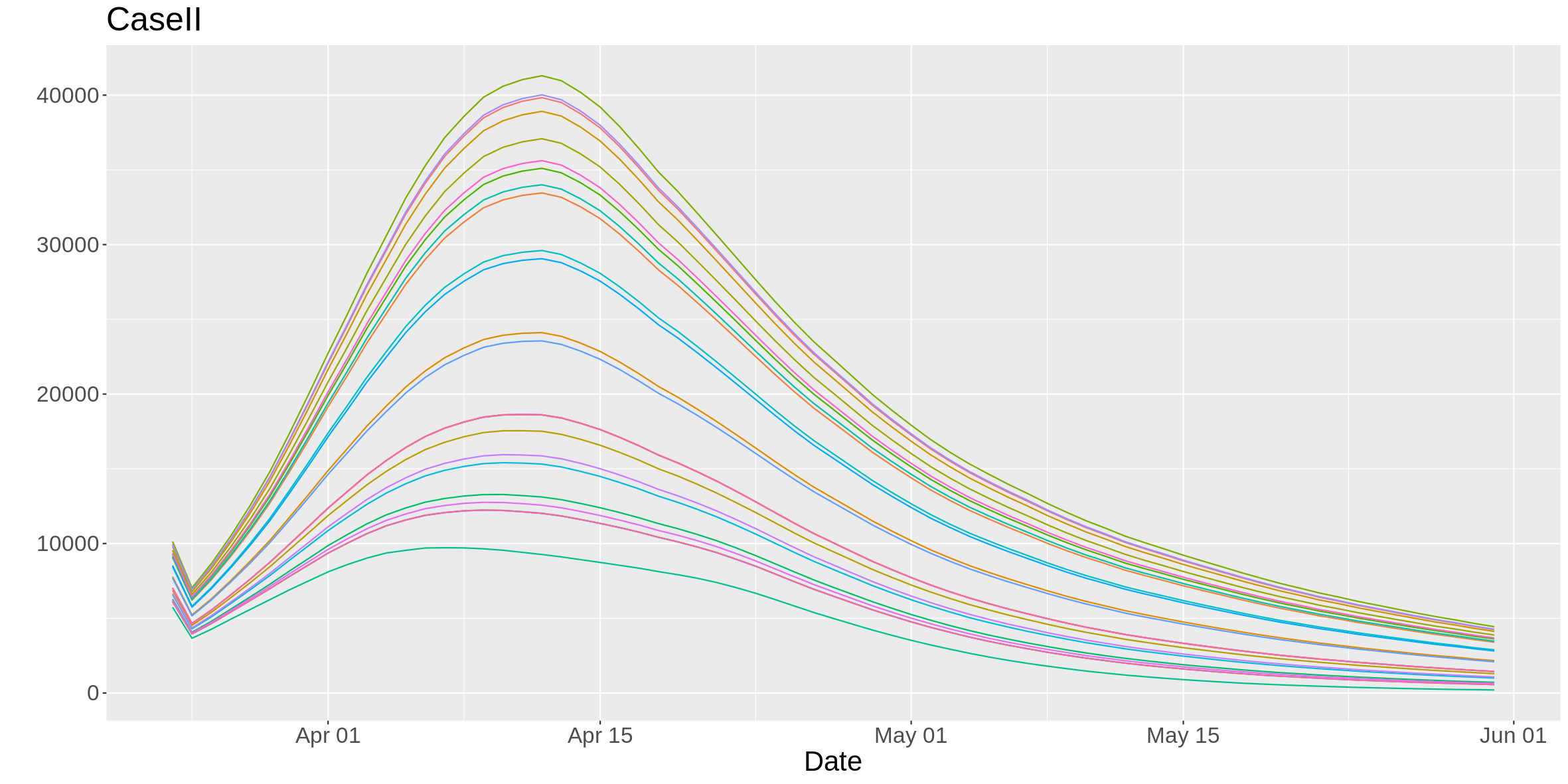}\\
		\includegraphics[width=65mm]{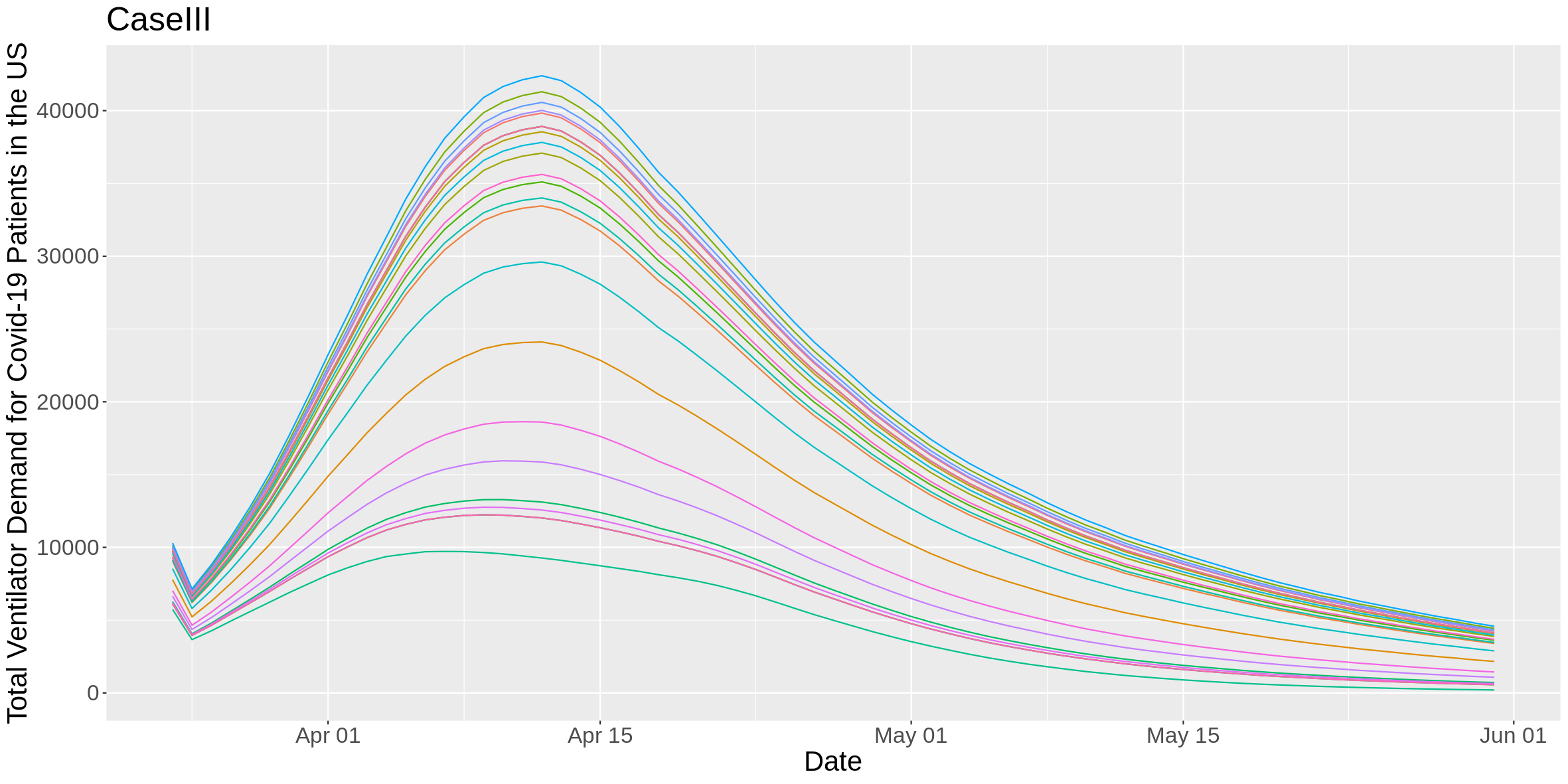} &   \includegraphics[width=65mm]{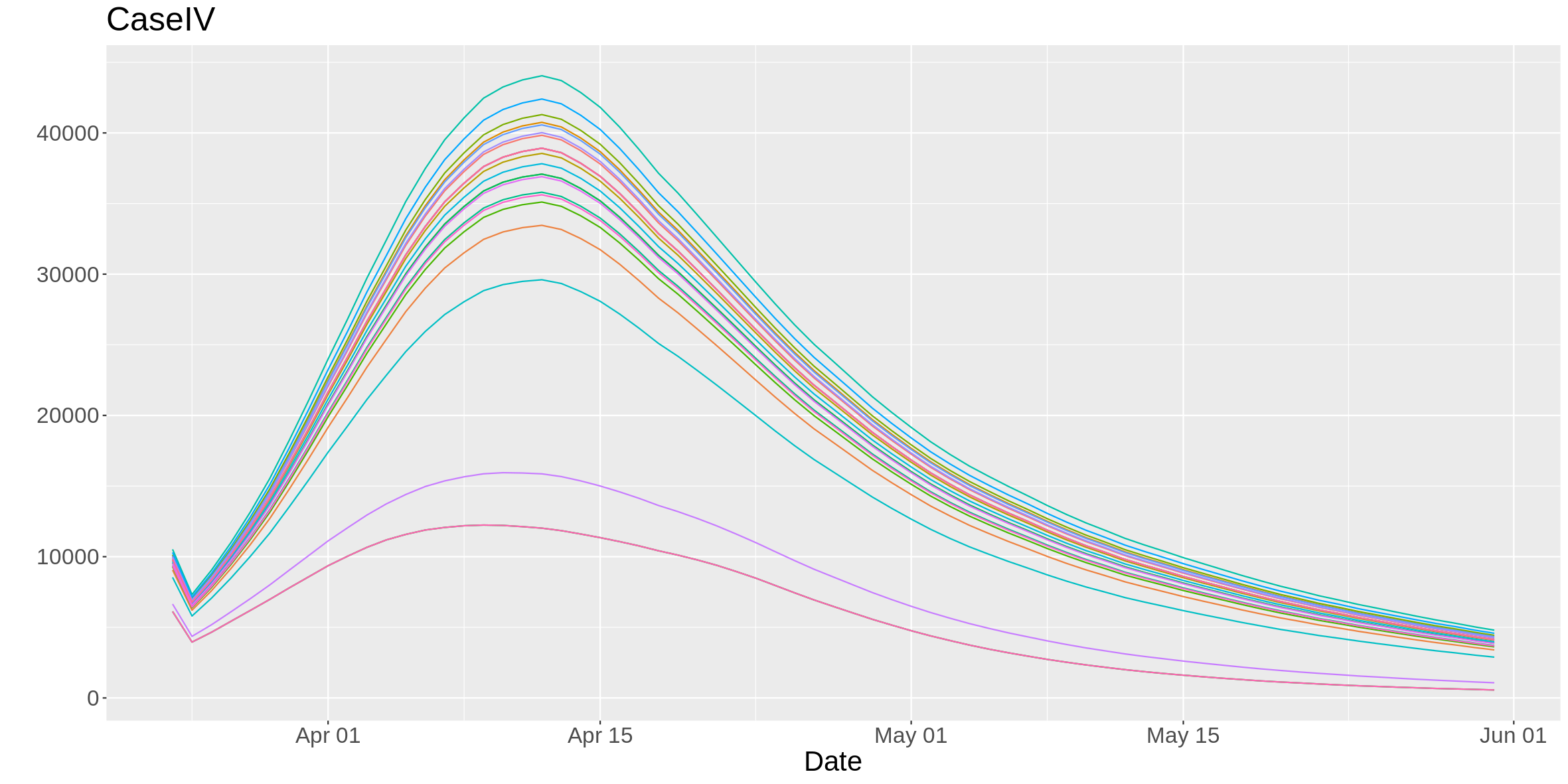} \\
		\multicolumn{2}{c}{(a) The US}\\
		\includegraphics[width=65mm]{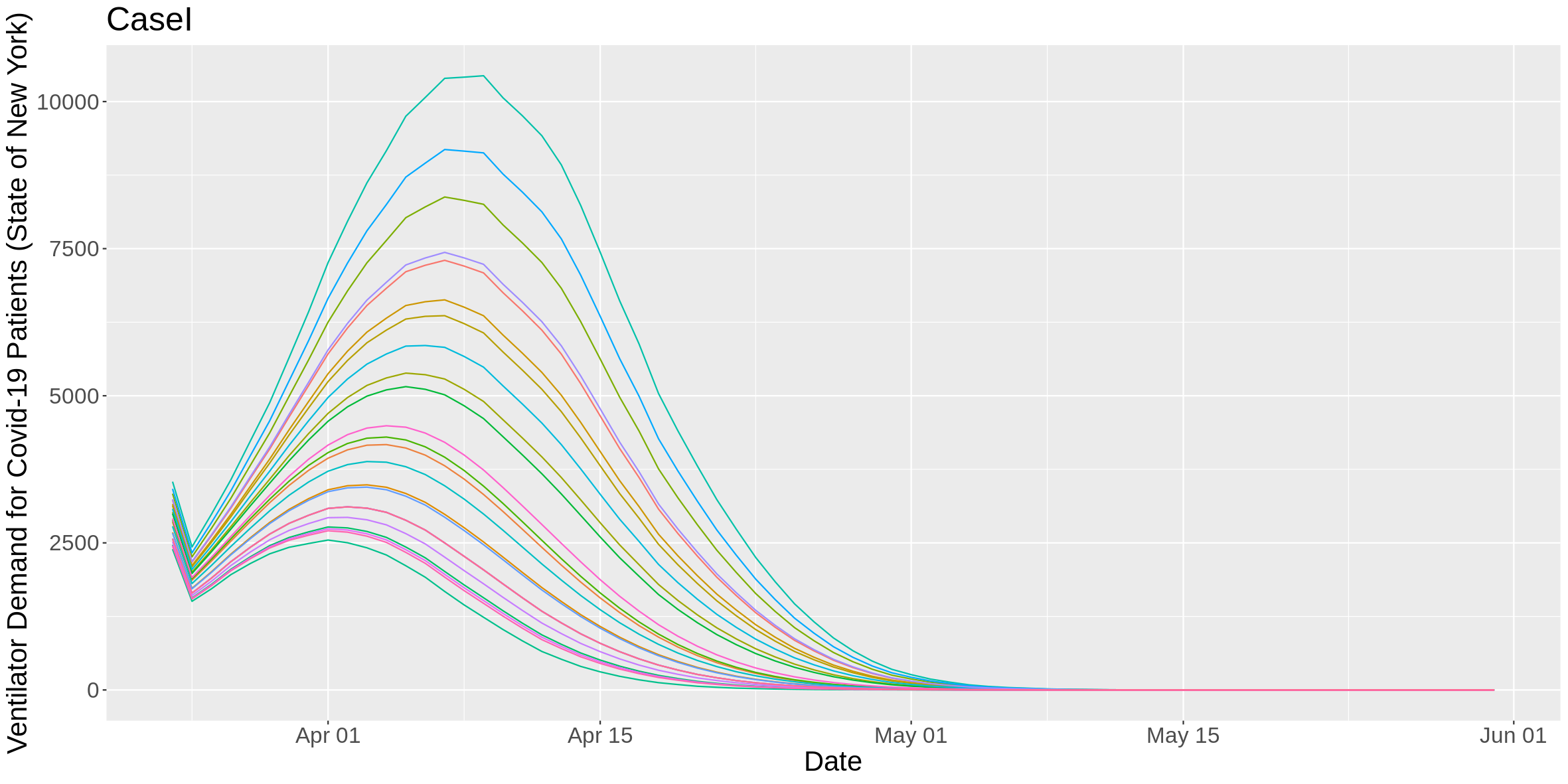} &   \includegraphics[width=65mm]{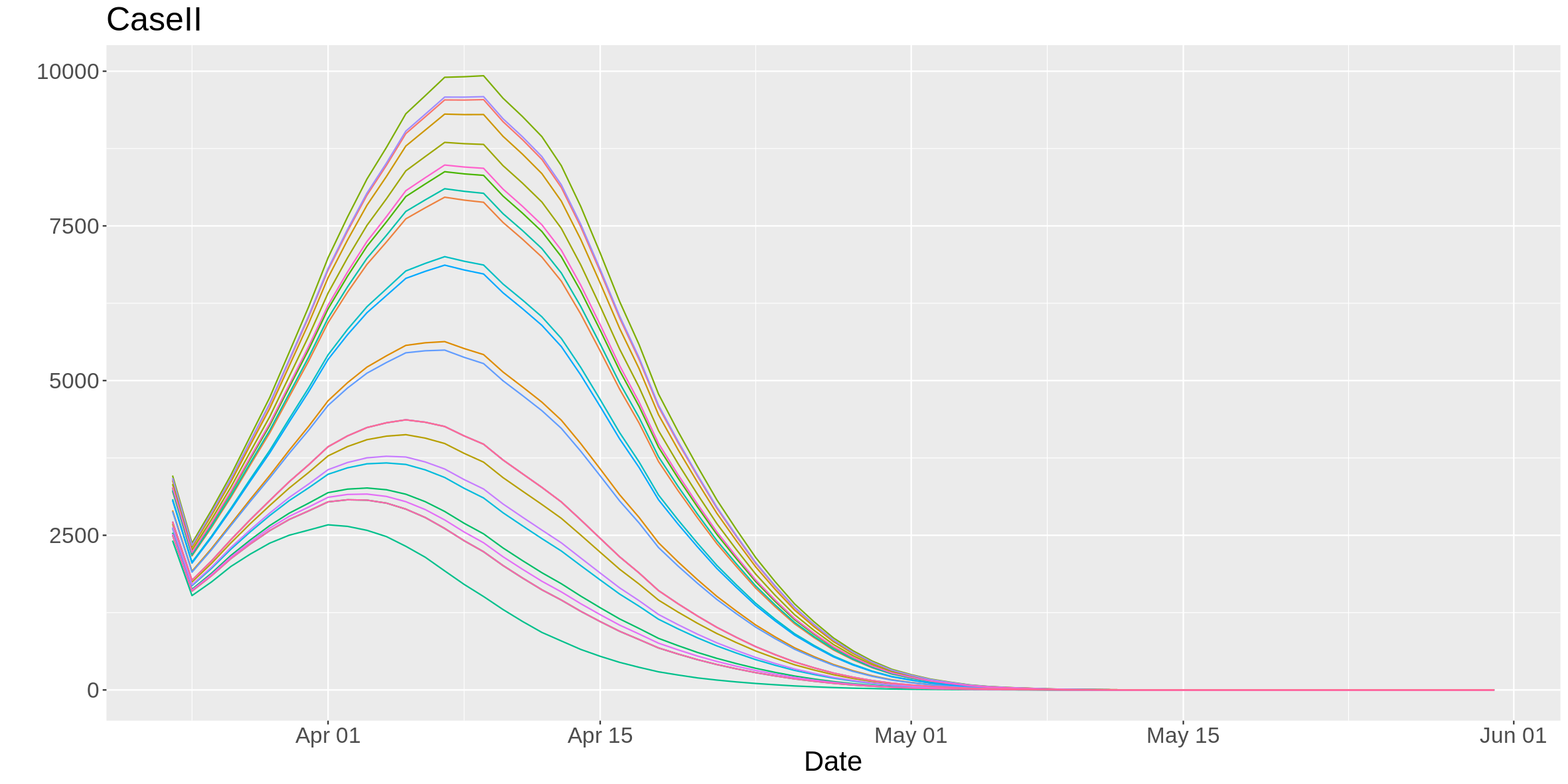}\\
		\includegraphics[width=65mm]{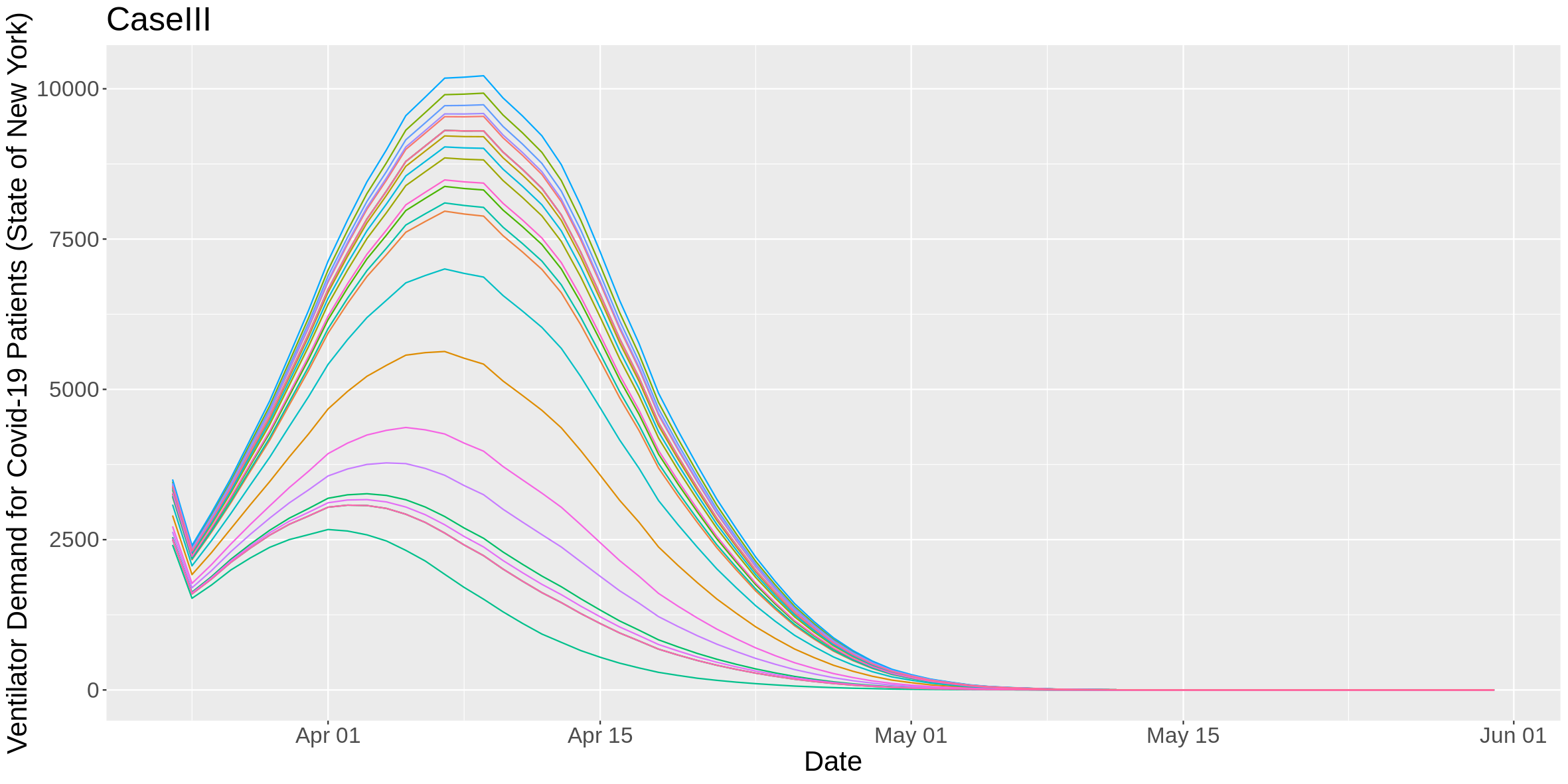} &   \includegraphics[width=65mm]{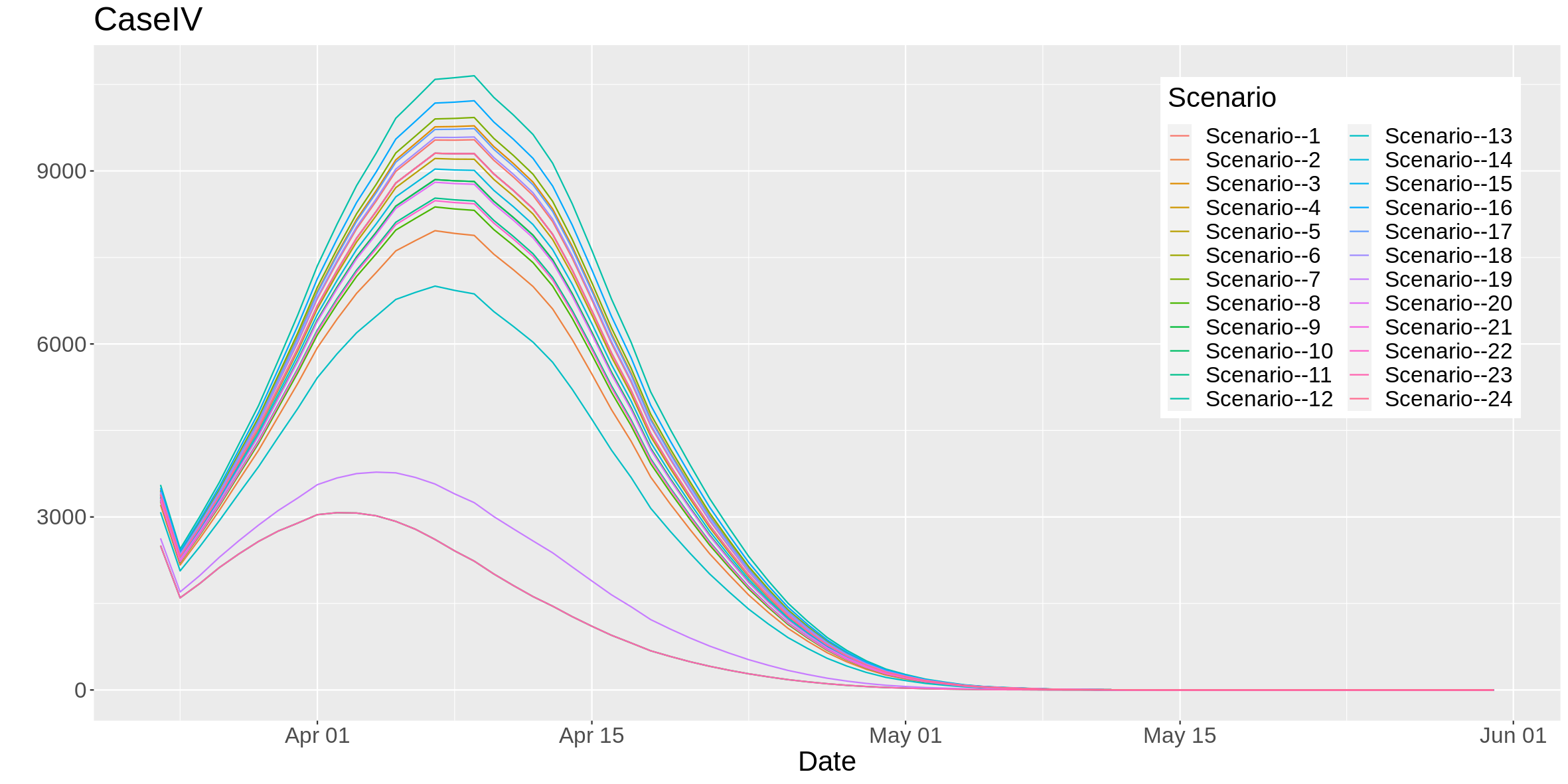}\\
		\multicolumn{2}{c}{(b) State of New York}\\
		\includegraphics[width=65mm]{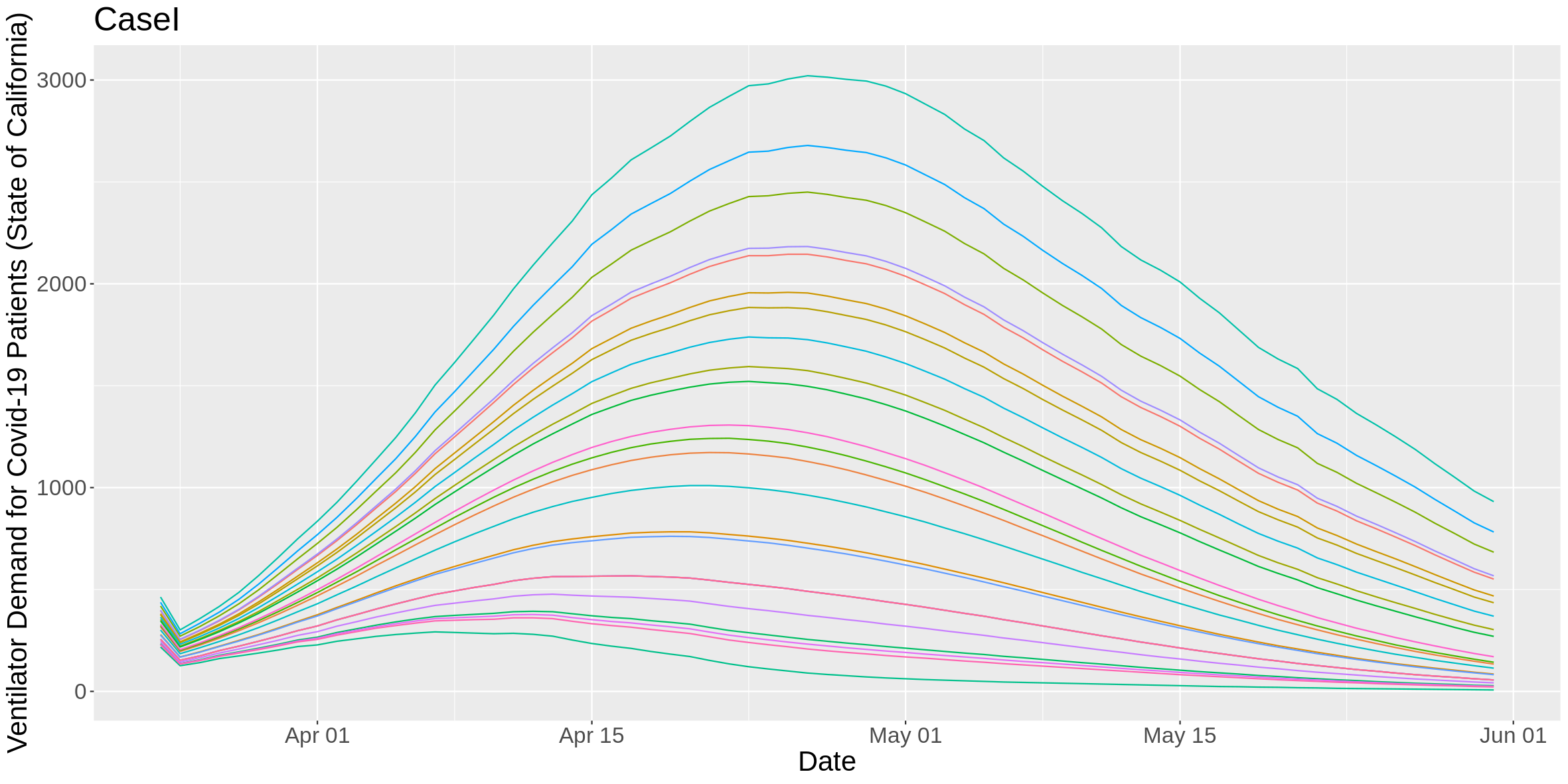} &   \includegraphics[width=65mm]{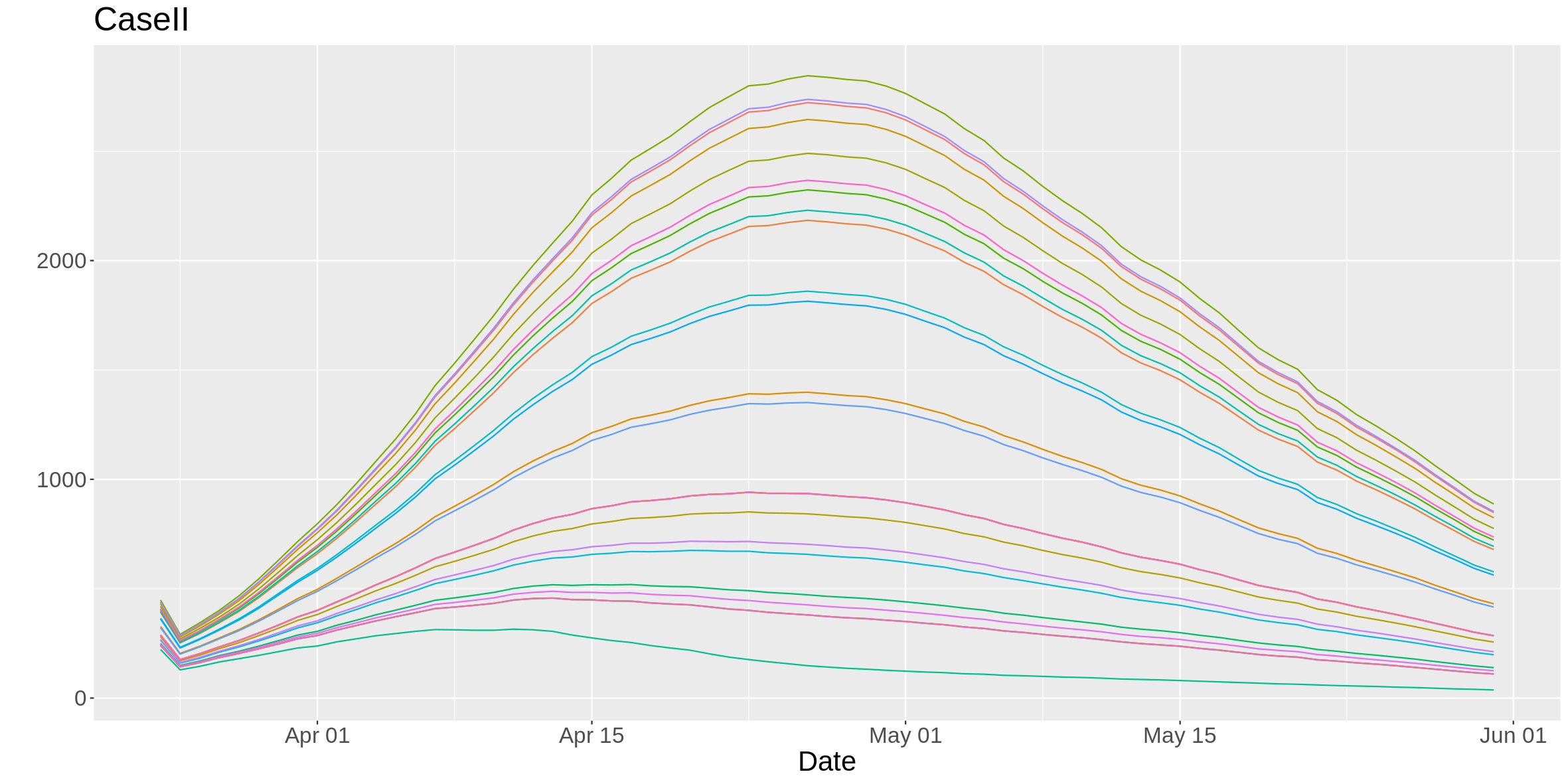}\\
		\includegraphics[width=65mm]{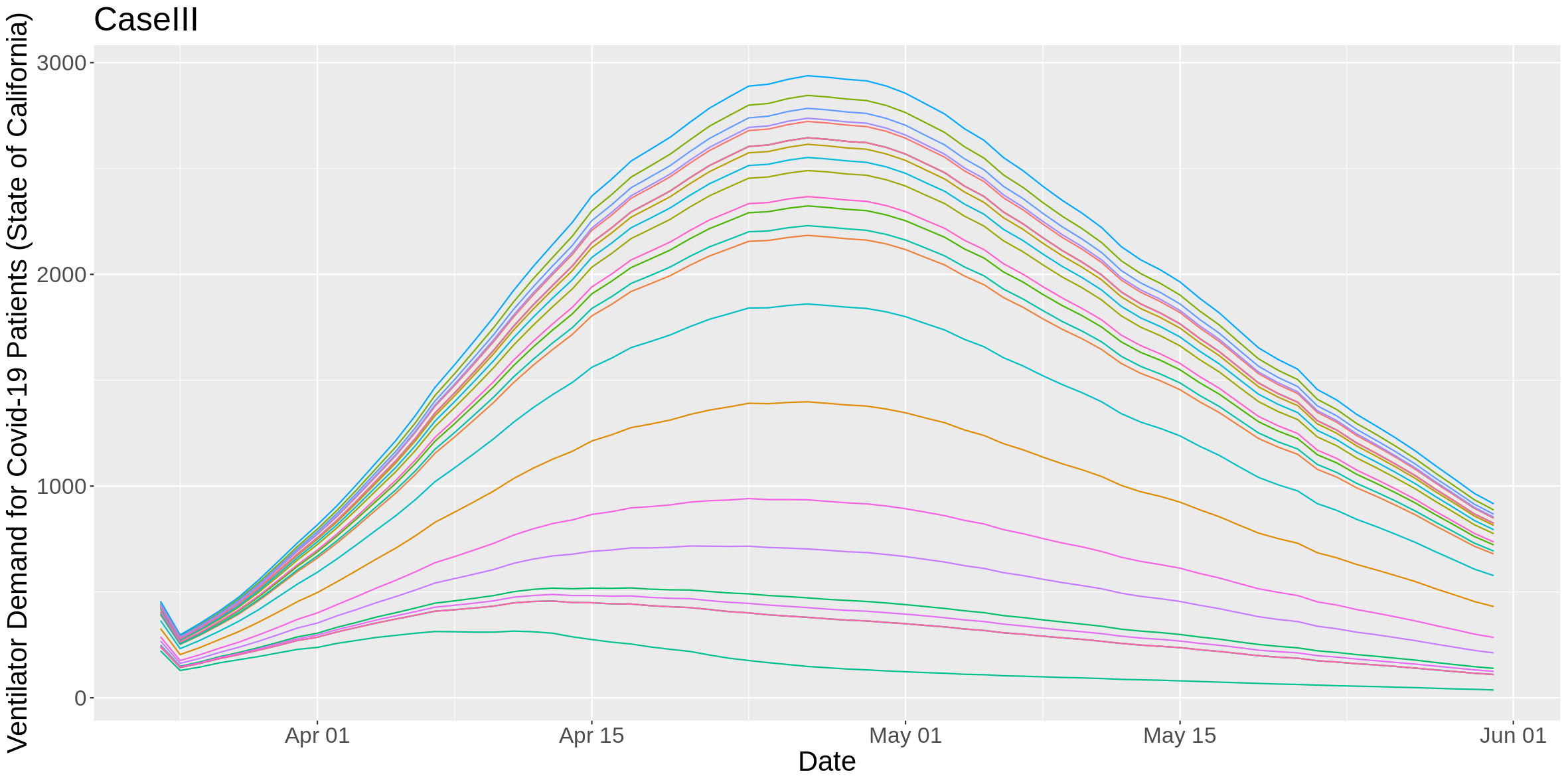} &   \includegraphics[width=65mm]{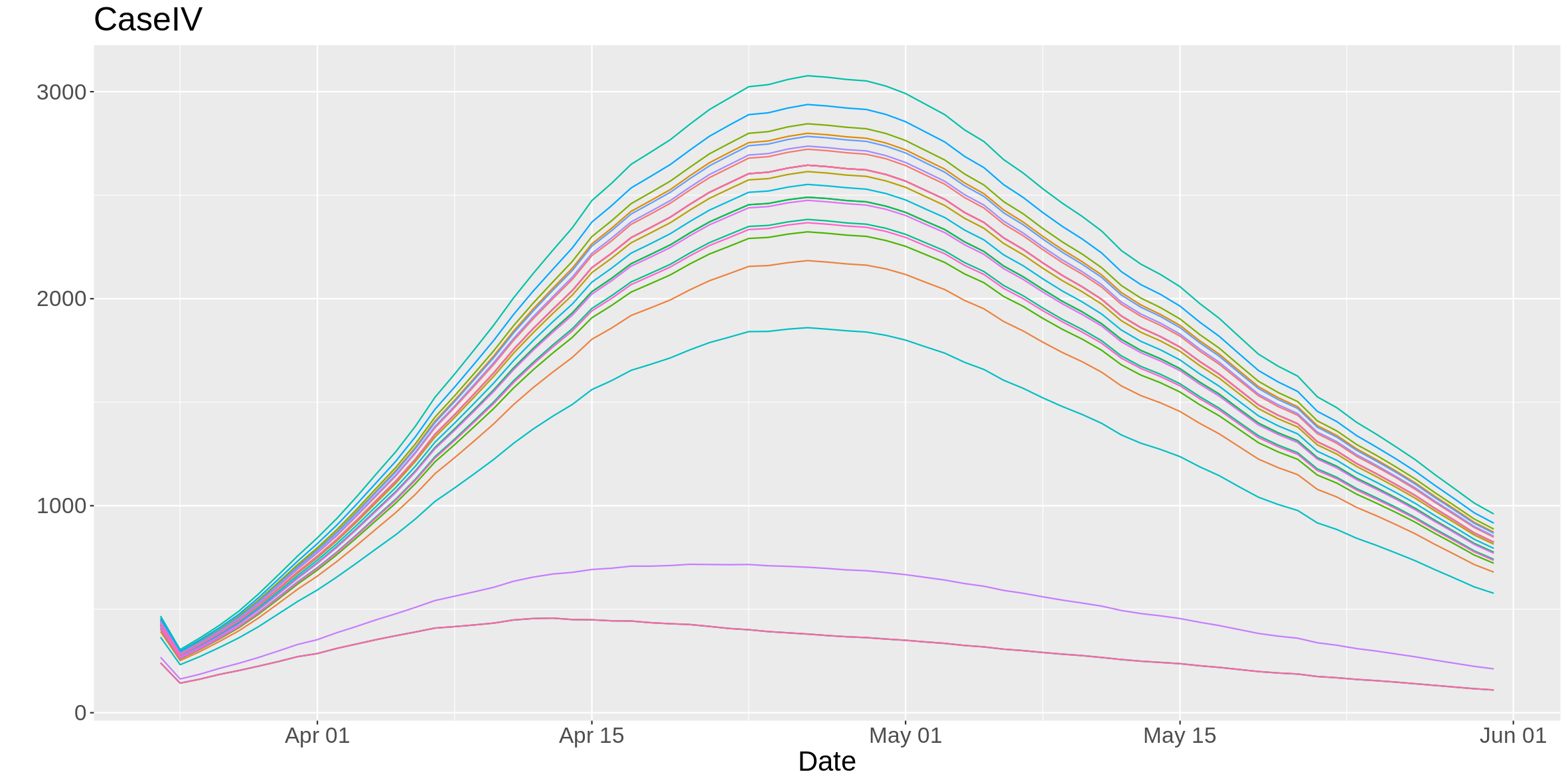} \\
		\multicolumn{2}{c}{(c) State of California}\\
	\end{tabular}
	\caption{Summary of generated scenarios (Cases I-IV) for the US and the States of New York and California, adapted from the data provided in~\cite{IHME2020}}
	\label{fig:1-US-NY-CA}
	% \includegraphics[width=130mm]{plots3/Demand.png}
	% \end{minipage}
	% \begin{minipage}{.3\linewidth}
	% \begin{tabular}{c}
	%       \
	% \end{tabular}
	% \end{minipage}
\end{figure}

For each of the four cases above, we generated 24 scenarios. Note that in each case, different quantities for the random demand  $\tilde{d}_{n,t}^{\omega}$,  $t \in \mtc{T}$, $n \in \mtc{N}$, and  $\omega \in \Omega$, might be generated. An illustration of the trajectory of demand scenarios over time is given in Figure~\ref{fig:1-US-NY-CA} for the US and the States of New York and  California.

% \begin{figure}%[h]
% \footnotesize
% \centering
% \includegraphics[width=0.95\textwidth]{plots2/US.png}
% \caption{Summary of generated scenarios (Case I) for the total ventilator demand in the U.S, adapted from the data provided in~\cite{IHME2020}.}
% \label{fig:1-US-SAA}
% \end{figure}

\subsection{Ventilator Inventory, Stockpile and Production}
Another key input to the planning model is the initial ventilator inventory. As of March 23, 2020, before the rapid rise of COVID-19 cases in NY, FEMA had about 20,000 ventilators in reserve, i.e., $I=20,000$. We used this for our model which suggests ventilator allocation decisions from 03/23/2020.

Estimates for the initial inventory of ventilators at different states were obtained from \cite{CovidCare}. These estimates are based on a hospital survey \cite{rubinson_vaughn_nelson_giordano_kallstrom_buckley_burney_hupert_mutter_handrigan_etal._2010,SCCM2020-Blog}. The estimates for new ventilator production were obtained  based on information provided at the US presidential briefings on 03/27/20 \cite{Trump032720}. These estimates suggest that the normal yearly ventilator production capacity is about 30,000 ventilators/year. However, under the US Defense Production Act, with the participation of additional companies, production of approximately 10,000 ventilators/month could be possible  \cite{Trump032720}. Using this information, for the baseline case we assumed that the current daily ventilator production rate is 100 ventilators/day; and it will be increased to 300 ventilators/day starting on April 15th.  

\subsection{Inventory Sharing Parameters}
Recall that in the model, parameter $\gamma$ is used to indicate the fraction of ventilators used  to care for non-COVID-19 patients. Additionally, a parameter $\tau$ is used in the model to estimate a state's willingness to share the fraction of their initial  COVID-19-use ventilators. Similarly, the parameter $\rho$ is used to control the state's risk-aversion to sending their idle ventilators to FEMA for use in a different state. 
We suppose that for all states $n$, $n \in \mtc{N}$, we have $\gamma_{n}=\gamma$, $\rho_{n}=\rho$, and $\tau_{n}=\tau$. 
In order to systematically study the ventilator allocations and shortfalls, we fixed the value of $\rho$ to $1.5$, and we used the following parameters: $\gamma \in \{50\%, 60\%, 75\%\}$ and $\tau \in \{0\%, 10\%, 25\%\}$.  

\subsection{Numerical Results}
For each setting $(\gamma, \tau)$, we solved model \eqref{opt:allocation_linear} under Cases I--IV. 
A summary of results is reported in Tables~\ref{T: summary} and \ref{T: summary_in_out}. 
We briefly describe the columns in these tables. 
Column ``Total" in Table~\ref{T: summary} denotes the total shortage, and is  calculated as 
\begin{equation*}
    \textrm{Total}:= 
\sum_{\omega \in \Omega} p^{\omega} \left[\sum_{t \in \mtc{T}}  \sum_{n\in \mtc{N}}   e_{n,t}^{\omega} \right]. 
\end{equation*}
Quantity ``Worst day" in column ``Worst day  $(t)$" denotes the shortage in the worst day, and is calculated as 
\begin{equation*}
    \textrm{Worst day}:= 
\max_{t \in \mtc{T}} \sum_{\omega \in \Omega} p^{\omega} \left[ \sum_{n\in \mtc{N}}   e_{n,t}^{\omega} \right], 
\end{equation*}
where $t$ denotes a day that the worst shortage happens, i.e., $t \in \argmax_{t \in \mtc{T}} \sum_{\omega \in \Omega} p^{\omega} \left[ \sum_{n\in \mtc{N}}   e_{n,t}^{\omega} \right]$. 
Moreover, quantity ``Worst day-state" in column ``Worst day-state  $(t)$" denotes the shortage in the worst day and state, and is calculated as 
\begin{equation*}
    \textrm{Worst day-state}:= 
\max_{t \in \mtc{T}} \max_{n\in \mtc{N}} \sum_{\omega \in \Omega} p^{\omega} e_{n,t}^{\omega}, 
\end{equation*}
where $(t,n) \in \argmax_{t \in \mtc{T}} \argmax_{n\in \mtc{N}} \sum_{\omega \in \Omega} p^{\omega} e_{n,t}^{\omega}$.

\begin{table}[!hb] 
\centering
\footnotesize
\begin{threeparttable} 
\caption{Ventilators' shortage summary under Cases I--IV.} 
\label{T: summary} 
\begin{tabular}{|l|l|l|l|l|l}
 \toprule 
$(\gamma, \tau)$ & Case & Total & Worst day ($t$) & Worst day-state  ($t$, $n$)   \\ 
 \midrule 
$(50\%, 0\%)$ &I & 112.21 & 22.12  (04/08/2020) &  15.38  (04/09/2020, New York)  \\ 
  &II & 0.00 & 0.00  &  0.00  \\ 
  &III & 23.29 & 8.88  (04/12/2020) &  3.12  (04/12/2020), Missouri)  \\ 
  &IV & 302.75 & 61.35  (04/12/2020) &  25.40  (04/09/2020, New York)  \\ 
 \midrule 
$(50\%, 25\%)$ &I & 0.00 & 0.00  &  0.00  \\ 
  &II & 0.00 & 0.00  &  0.00  \\ 
  &III\tnote{\textdagger}  & 3.58 & 2.04  (04/10/2020) &  1.96  (04/10/2020, New Jersey)  \\ 
  &IV & 0.00 & 0.00  &  0.00  \\ 
 \midrule 
$(50\%, 50\%)$ &I & 0.00 & 0.00  &  0.00  \\ 
  &II & 0.00 & 0.00  &  0.00  \\ 
  &III & 0.00 & 0.00  &  0.00  \\ 
  &IV & 0.00 & 0.00  &  0.00  \\ 
 \midrule 
$(60\%, 0\%)$ &I & 611.21 & 75.37  (04/10/2020) &  41.67  (04/09/2020, New York)  \\ 
  &II & 95.95 & 22.15  (04/12/2020) &  7.30  (04/09/2020, New York)  \\ 
  &III & 2525.46 & 288.38  (2020-04-09) &  183.50  (04/09/2020, New York)  \\ 
  &IV\tnote{\textdagger}  & 2000.80 & 296.75  (04/12/2020) &  139.70  (04/09/2020, New York)  \\ 
 \midrule 
$(60\%, 25\%)$ &I\tnote{\textdagger}  & 16.42 & 3.67  (04/12/2020) &  2.38  (04/08/2020, Michigan)  \\ 
  &II & 0.00 & 0.00  &  0.00  \\ 
  &III & 0.00 & 0.00  &  0.00  \\ 
  &IV\tnote{\textdagger}  & 2072.70 & 253.30  (04/12/2020) &  92.15  (04/08/2020, Michigan)  \\ 
 \midrule 
$(60\%, 50\%)$ &I & 0.00 & 0.00  &  0.00  \\ 
  &II & 0.00 & 0.00  &  0.00  \\ 
  &III\tnote{\textdagger}  & 157.00 & 24.79  (04/12/2020) &  8.08  (04/07/2020, Michigan)  \\ 
  &IV\tnote{\textdagger}  & 32.10 & 6.10  (04/12/2020) &  4.30  (04/08/2020, Michigan)  \\ 
 \midrule 
$(75\%, 0\%)$ &I\tnote{\textdagger}  & 4395.46 & 428.08  (04/12/2020) &  153.92  (04/07/2020, New York)  \\ 
  &II\tnote{\textdagger}  & 2877.37 & 299.10  (04/12/2020) &  119.83  (04/07/2020, New York)  \\ 
  &III\tnote{\textdagger}  & 15748.62 & 1548.04  (04/12/2020) &  642.54  (04/07/2020, New York)  \\ 
  &IV\tnote{\textdagger}  & 28529.72 & 2693.77  (04/12/2020) &  1237.10  (04/07/2020, New York)  \\ 
 \midrule 
$(75\%, 25\%)$ &I & 4260.38 & 372.21  (04/12/2020) &  169.96  (04/07/2020, New York)  \\ 
  &II\tnote{\textdagger}  & 3197.37 & 305.92  (04/12/2020) &  123.72  (04/07/2020, New York)  \\ 
  &III\tnote{\textdagger}  & 15026.58 & 1368.63  (04/09/2020) &  646.88  (04/07/2020, New York)  \\ 
  &IV\tnote{\textdagger}  & 26990.43 & 2436.63  (04/12/2020) &  1168.08  (04/07/2020, New York)  \\ 
 \midrule 
$(75\%, 50\%)$ &I\tnote{\textdagger}  & 3667.96 & 339.42  (04/12/2020) &  138.46  (04/07/2020, New York)  \\ 
  &II\tnote{\textdagger}  & 2336.97 & 217.83  (04/09/2020) &  111.38  (04/07/2020, New York)  \\ 
  &III\tnote{\textdagger}  & 13300.00 & 1244.92  (04/09/2020) &  591.54  (04/07/2020, New York)  \\ 
  &IV\tnote{\textdagger}  & 24828.15 & 2264.55  (04/09/2020) &  1081.03  (04/07/2020, New York)  \\ 
\bottomrule 
\end{tabular} 
\begin{tablenotes}\footnotesize 
\item[\textdagger] Reached the one-hour time limit. The reported results correspond to the best integer solution found.
\end{tablenotes} 
\end{threeparttable} 
\end{table}

We also analyzed the ventilators' reallocation to/from different states for the setting $(\gamma, \tau)=(0.75, 0)$, which is the most dramatic case we considered from the inventory and stockpile perspectives. We report a summary of results  in Table  \ref{T: summary_in_out} under the two worst situations, Cases III (Mildy Worse than Average) and IV (Severe).
Column ``Total inflow" in this table denotes the total incoming ventilators to a state $n \in \mtc{N}$ from FEMA, and is  calculated as 
\begin{equation*}
    \textrm{Total inflow}:= 
\sum_{t \in \mtc{T}}  \sum_{\omega \in \Omega} p^{\omega} x_{n,t}^{\omega}.
\end{equation*}
Similarly, 
column ``Total outflow" denotes the total outgoing ventilators from a state $n \in \mtc{N}$ to FEMA, and is  calculated as 
\begin{equation*}
    \textrm{Total outflow}:= 
\sum_{t \in \mtc{T}}  \sum_{\omega \in \Omega} p^{\omega} z_{n,t}^{\omega}.
\end{equation*}
Also, column ``Net flow" represents the difference between ``Total inflow" and ``Total outflow".

%\caption{Inflow and outflow of ventilators under Case III (Severe) and Case IV (Extreme), with $\gamma=0.75$, $\rho=1.5$, and $\tau=0$.} 

\begin{table} 
\caption{Inflow and outflow of ventilators under Case III (Mildy Worse than Average) and Case IV (Severe), with $(\gamma, \tau)=(0.75, 0)$.} 
\label{T: summary_in_out} 
\begin{adjustbox}{width=\textwidth} 
\begin{tabular}{lllllll}
 \toprule 
 & \multicolumn{3}{c}{Case III} & \multicolumn{3}{c}{Case IV} \\ 
 \cmidrule(lr){2-4}\cmidrule(lr){5-7}State  & Total inflow & Total outflow & Net flow & Total inflow & Total outflow & Net flow \\ 
 \midrule 
 Alabama  & 379.75 & 0.00 & 379.75  & 460.85 & 7.60 & 453.25  \\ 
  Alaska  & 382.62 & 0.54 & 382.08  & 295.22 & 0.00 & 295.22  \\ 
  Arizona  & 274.71 & 1.00 & 273.71  & 320.77 & 2.10 & 318.67  \\ 
  Arkansas  & 29.96 & 0.00 & 29.96  & 119.17 & 0.00 & 119.17  \\ 
  California  & 5755.87 & 360.62 & 5395.25  & 4343.08 & 304.20 & 4038.88  \\ 
  Colorado  & 240.87 & 0.00 & 240.87  & 397.42 & 0.00 & 397.42  \\ 
  Connecticut  & 514.83 & 0.00 & 514.83  & 526.57 & 84.20 & 442.37  \\ 
  Delaware  & 111.88 & 0.00 & 111.88  & 60.13 & 0.00 & 60.13  \\ 
  District of Columbia  & 171.25 & 0.00 & 171.25  & 179.23 & 0.00 & 179.23  \\ 
  Florida  & 72.58 & 0.00 & 72.58  & 21.15 & 0.00 & 21.15  \\ 
  Georgia  & 624.08 & 1.54 & 622.54  & 837.52 & 1.80 & 835.72  \\ 
  Hawaii  & 32.50 & 0.00 & 32.50  & 253.68 & 0.00 & 253.68  \\ 
  Idaho  & 178.42 & 0.00 & 178.42  & 159.13 & 1.00 & 158.13  \\ 
  Illinois  & 547.58 & 30.42 & 517.17  & 828.43 & 7.97 & 820.47  \\ 
  Indiana  & 1306.50 & 88.08 & 1218.42  & 1885.82 & 88.00 & 1797.82  \\ 
  Iowa  & 62.37 & 0.00 & 62.37  & 101.70 & 0.00 & 101.70  \\ 
  Kansas  & 78.67 & 0.12 & 78.54  & 163.38 & 0.00 & 163.38  \\ 
  Kentucky  & 353.17 & 0.00 & 353.17  & 155.55 & 0.00 & 155.55  \\ 
  King and Snohomish Counties, WA  & 1072.08 & 2.00 & 1070.08  & 1618.42 & 195.72 & 1422.70  \\ 
  Louisiana  & 1187.46 & 14.92 & 1172.54  & 1351.48 & 32.80 & 1318.68  \\ 
  Maine  & 139.58 & 7.50 & 132.08  & 128.20 & 0.10 & 128.10  \\ 
  Maryland  & 410.75 & 0.00 & 410.75  & 327.03 & 3.35 & 323.68  \\ 
  Massachusetts  & 1304.42 & 42.08 & 1262.33  & 1675.13 & 55.05 & 1620.08  \\ 
  Michigan  & 3003.46 & 365.92 & 2637.54  & 3181.92 & 462.43 & 2719.48  \\ 
  Minnesota  & 347.08 & 0.33 & 346.75  & 172.70 & 0.00 & 172.70  \\ 
  Mississippi  & 38.29 & 0.00 & 38.29  & 22.85 & 0.00 & 22.85  \\ 
  Missouri  & 2235.88 & 558.04 & 1677.83  & 2430.73 & 332.15 & 2098.58  \\ 
  Montana  & 98.62 & 0.00 & 98.62  & 72.12 & 0.00 & 72.12  \\ 
  Nebraska  & 6.67 & 0.00 & 6.67  & 19.72 & 0.00 & 19.72  \\ 
  Nevada  & 630.37 & 13.08 & 617.29  & 604.22 & 9.20 & 595.02  \\ 
  New Hampshire  & 95.21 & 0.04 & 95.17  & 118.00 & 0.53 & 117.47  \\ 
  New Jersey  & 2087.25 & 259.42 & 1827.83  & 2650.07 & 385.95 & 2264.12  \\ 
  New Mexico  & 69.50 & 0.00 & 69.50  & 53.87 & 0.30 & 53.57  \\ 
  New York  & 9140.00 & 2086.33 & 7053.67  & 7769.53 & 420.08 & 7349.45  \\ 
  North Carolina  & 376.58 & 3.54 & 373.04  & 358.77 & 0.00 & 358.77  \\ 
  North Dakota  & 467.04 & 60.71 & 406.33  & 64.17 & 0.00 & 64.17  \\ 
  Ohio  & 66.54 & 0.00 & 66.54  & 133.48 & 0.00 & 133.48  \\ 
  Oklahoma  & 99.54 & 0.00 & 99.54  & 80.67 & 5.85 & 74.82  \\ 
  Oregon  & 66.38 & 0.00 & 66.38  & 116.92 & 0.00 & 116.92  \\ 
  Other Counties, WA  & 1050.92 & 46.54 & 1004.38  & 1151.35 & 44.85 & 1106.50  \\ 
  Pennsylvania  & 138.21 & 0.00 & 138.21  & 149.15 & 6.65 & 142.50  \\ 
  Rhode Island  & 53.29 & 0.00 & 53.29  & 73.37 & 0.00 & 73.37  \\ 
  South Carolina  & 45.29 & 0.00 & 45.29  & 100.35 & 0.00 & 100.35  \\ 
  South Dakota  & 26.67 & 0.00 & 26.67  & 79.87 & 0.05 & 79.82  \\ 
  Tennessee  & 1058.75 & 507.54 & 551.21  & 480.03 & 0.00 & 480.03  \\ 
  Texas  & 768.62 & 0.33 & 768.29  & 331.80 & 0.00 & 331.80  \\ 
  Utah  & 74.33 & 0.00 & 74.33  & 94.27 & 0.00 & 94.27  \\ 
  Vermont  & 732.62 & 85.88 & 646.75  & 280.43 & 93.10 & 187.33  \\ 
  Virginia  & 2750.17 & 433.25 & 2316.92  & 1820.38 & 37.10 & 1783.28  \\ 
  Washington  & 144.71 & 0.00 & 144.71  & 216.08 & 0.00 & 216.08  \\ 
  West Virginia  & 6.17 & 1.46 & 4.71  & 3.95 & 0.00 & 3.95  \\ 
  Wisconsin  & 483.50 & 0.00 & 483.50  & 48.90 & 0.00 & 48.90  \\ 
  Wyoming  & 51.17 & 0.00 & 51.17  & 56.65 & 0.00 & 56.65  \\ 
 \bottomrule 
\end{tabular} 
\end{adjustbox} 
\end{table}

\subsection{Discussion}
The results in Table~\ref{T: summary} suggest that when up to 50\% of a state's ventilator inventory is used for non-COVID-19 patients, FEMA's current stockpile of 20,000 ventilators is sufficient  to meet the demand imposed by COVID-19 patients. This ventilator use threshold increases to 60\% of non-COVID-19 patients,  if states are willing to share up to 50\% of their excess inventory with other states. However, if no such sharing is considered, then a moderate number of  ventilators (approximately 300) will be required beyond FEMA's current stockpile to meet demand in Cases I-IV. 

The ventilator availability situation gets worse in the case where 75\% (or greater \%) of the available ventilators must be used for non-COVID-19 patients. In this scenario, in Case III (Mildy Worse than Average) and Case IV (Severe) the inventory shortfall on the worst day (04/12/2020) is between 1,500-2,700. This shortfall decreases moderately to 1,250-2,250 if states are willing to share part of their initial ventilator inventory. 

%\begin{enumerate}
%    \item[Conclusion 1.] As the states are willing to share a portion of their original stocks, ??? lives can be saved.
%    \item[Conclusion 2.] As the risk aversion of the states reduce, and they are willing to share the ventilators, ??? lives can be saved.
%\end{enumerate}

\section{Concluding Remarks}
We have presented a model for procuring and sharing life-saving resources whose demand is stochastic. The demand arising from different entities (states) peaks at different times, and it is important to meet as much of this demand as possible to save lives. Each participating state is risk averse to sharing their excess inventory at any given time, and this risk-aversion is captured by using a safety threshold parameter. Specifically, the developed model is applicable to the current COVID-19 pandemic, where many US states are in dire need of mechanical ventilators to provide life-support to severely- and critically-ill patients. Computations were performed using realistic ventilator need forecasts and availability under a wide combination of parameter settings. 

Our findings suggest that the fraction of currently available ventilators that are to be used for non-COVID-19 patients strongly impacts state and national ability to meet demand arising from COVID-19 patients. When more than 40\% of the existing inventory is available for COVID-19 patients, the national stockpile is sufficient to meet the demand. %\Add what sharing conditions this is under. Eg- assuming states are willing to share __% of their excess inventory at any given time. 
However, if less than 25\% of the existing inventory is available for COVID-19 patients, the current national stockpile and the anticipated production may not be sufficient under extreme demand scenarios. As expected, the magnitude of this shortfall increases when one considers more and more extreme demand scenarios.  
%\ or greater risk aversion to sharing on the part of the states? Add shortfall that will occur under more extreme demand scenarios if the states are not willing to share.  

Overall, the model developed in this paper can be used as a planning tool/framework by state and federal agencies in acquiring and allocating ventilators to meet national demand. The results reported in this paper can also provide a guide to states in planning for their ventilator needs. We, however, emphasize that these results are based on certain modeling assumptions. This include the process of demand forecast scenario generation, estimates of initial ventilator inventory, and future production quantities. Each one of these, as well as other model parameters, can be changed in the model input to obtain more refined results. Nevertheless, an important finding is that a state's willingness to share its idle inventory can help address overall shortfall. 

While this paper has focused on ventilator needs in the US, such a model can also be adapted for use in international supply-chain coordination of equipment such as ventilators across countries. COVID-19 is expected to have different peak dates and demand cycles in other counties, and one or two additional disease spread cycles are likely until an effective vaccine becomes available.  

In conclusion, we point out that the model developed in this paper has a one-time planning decision, i.e., there are no ``wait-and-see"  decisions in the model. One can also formulate the ventilator allocation problem  as a time-dynamic multistage stochastic program, where the decision maker can make recourse decisions as time evolves based on the information available so far on the stochastic demands and past decisions. We are currently working on such an extension.

\newpage
\bibliographystyle{abbrvnat}
\bibliography{main}

\end{document}